\documentclass[a4paper, 11pt]{article} 
\usepackage{amsfonts,amsmath,amssymb,amscd}
\usepackage{latexsym}
\usepackage{graphicx}
\usepackage{color} 
 \usepackage{hyperref}

\newtheorem{Theo}{Theorem}[section]
\newtheorem{Coro}[Theo]{Corollary}
\newtheorem{Lemma}[Theo]{Lemma}

\newtheorem{Propo}[Theo]{Proposition}

\newtheorem{Defi}[Theo]{Definition}
\newtheorem{Rema}[Theo]{Remark}
\newtheorem{Exam}[Theo]{Example}
\numberwithin{equation}{section}

\newcommand{\beqr}{\begin{eqnarray}}
\newcommand{\eeqr}{\end{eqnarray}}

\newcommand{\beq}{\begin{eqnarray*}}
\newcommand{\eeq}{\end{eqnarray*}}

\newcommand{\bq}{\begin{equation}}
\newcommand{\eq}{\end{equation}}

\newenvironment{preuve}[1][]
{\vskip 2mm  {\it \bf Proof#1. }}{$\Box$ \vskip 2mm}

\newcommand{\bpr}{\begin{preuve}}
\newcommand{\epr}{\end{preuve}}

\newcommand{\bef}{\mathcal{F}}
\newcommand{\be}{\mathcal{E}}

\newcommand{\bh}{\mathcal{H}}
\newcommand{\bi}{\mathcal{I}}

\newcommand{\Nn}{\mathbb{N}}
\newcommand{\R}{\mathbb{R}}

\newcommand{\cg}{\langle}
\newcommand{\cd}{\rangle}

\newcommand{\ep}{\epsilon}
\newcommand{\si}{\sigma}
\newcommand{\Si}{\Sigma}

\newcommand{\udn}{\frac{1}{2^n}}
\newcommand{\ud}{\frac{1}{2}}

\newcommand{\upn}{\frac{1}{(2\pi)^n}}

\newcommand{\ck}{\chi_K}
\newcommand{\ldn}{L^2(\R^n)}
\newcommand{\lkn}{L^2_K(\R^n)}

\newcommand{\sumun}{\sum_{i=1}^{n}}
\newcommand{\un}{\{1, \cdots, n\}}
\newcommand{\zn}{\{0, \cdots, n-1\}}
\newcommand{\izn}{i\in \{0, \cdots, n-1\}}

\newcommand{\lumm}{L^{-\frac{1}{m}  }}

\newcommand{\lnmm}{\frac{1}{L^\frac{n}{m}  }}
\newcommand{\lndmm}{L^{-\frac{n}{2m}  }}
\newcommand{\lndm}{L^{\frac{n}{2m}  }}

\newcommand{\gme}{\Gamma(M,E)}

\newcommand{\sz}{s^{-1}(0)}

\newcommand{\syp}{\sigma_P}
\newcommand{\equiL}{\underset{L\to +\infty}{\sim}}

\newcommand{\toL}{\underset{L\to +\infty}{\to}}
\newcommand{\linf}{L\to +\infty}
\newcommand{\dl}{\Delta_L}
\newcommand{\Linf}{L^{\infty}}
\newcommand{\ul}{\mathbb U_L}
\newcommand{\uld}{\ul\setminus \Delta_L}
\newcommand{\met}{Met_{|dy|}(M)}
\newcommand{\rn}{\R^n}

\newcommand{\kx}{K_x}
\newcommand{\rks}{p^K_\Si}
\newcommand{\pxs}{p^x_\Si}

\newcommand{\tx}{T_x M}
\newcommand{\iso}{Isom_g (\rn, \tx)}
\newcommand{\csp}{c_\Sigma(P)}
\newcommand{\brl}{B_g(x,RL^{-\frac{1}{m}})}

\newcommand{\beu}{\begin{enumerate}}
\newcommand{\eeu}{\end{enumerate}}
\newcommand{\twf}{\mathcal T_{(W,f)}     }
\newcommand{\bt}{\mathcal T}

\newcommand{\pif}{\frac{\partial f}{\partial x_i}}
\newcommand{\wf}{(W,f)}
\newcommand{\rwf}{R_{(W,f)}}
\newcommand{\tkwf}{\tau^K_{(W,f)}}
\newcommand{\xz}{x_0}
\newcommand{\phx}{\phi_{x_0}^{-1}}
\newcommand{\px}{\phi_{x_0}}
\newcommand{\wfs}{(W,f_\Si)}
\newcommand{\bzr}{B(0,R_{\wfs}\lumm)}
\newcommand{\bgr}{B_g(\xz,R_{\wfs}\lumm)}

\newcommand{\fr}{\frac}
\newcommand{\pedn}{\lfloor \frac{n}{2}+1 \rfloor}
\newcommand{\dn}{\frac{n}{2}}
\newcommand{\spn}{\sqrt{\pi}^n}
\newcommand{\sdpn}{\sqrt{2\pi}^n}
\newcommand{\Rpe}{\R^*_+}
\newcommand{\rkr}{\rho_K(R)}
\newcommand{\tkjr}{\theta^j_K (R)}
\newcommand{\cep}{\frac{c}{2\eta}}

\newcommand{\cik}{C^\infty_0(K)}
\newcommand{\srn}{S(\rn)}
\newcommand{\befm}{\bef^{-1}}
\newcommand{\skn}{S_K(\R^n)}

\begin{document}

\baselineskip=17pt


\title{Universal components of random nodal sets}

\author{Damien Gayet
\and 
Jean-Yves Welschinger}

\date{\today}

\maketitle


\begin{abstract}
We give, as $L$ grows to infinity, an explicit lower bound 
of order $L^{\frac{n}m}$
for the expected Betti numbers of the vanishing locus of a random
linear combination of eigenvectors 
of $P$ with eigenvalues below $L$. 
Here, $P$ denotes
an elliptic self-adjoint pseudo-differential operator
of order $m>0$, bounded from below and acting 
on the sections of a Riemannian line bundle over 
a smooth closed $n$-dimensional manifold  $M$ equipped with some Lebesgue measure.
In fact, for every
closed hypersurface $\Si$ of $\rn$, we prove that there exists a positive constant $p_\Si$ depending only on $\Si$, such that for every  large enough $L$ and every $x\in M$,  a component diffeomorphic to $\Si$ appears
with probability at least $p_\Si$ in the vanishing locus of a random section and
in the ball of radius $L^{-\frac{1}{m}}$ centered at $x$. 
These results apply in particular to  Laplace-Beltrami and  Dirichlet-to-Neumann operators.\\


\end{abstract}

\section*{Introduction}
Let $M$ be a  smooth closed manifold of positive dimension $n$ and $E$ be a real line bundle over $M$.
We equip $M$ with a Lebesgue measure $|dy|$, that is a positive measure that can be locally expressed
as the absolute value of some smooth volume form, and 
$E$ with a Riemannian metric 
$h_E$. These induce a $L^2$-scalar product on the space $\Gamma(M,E)$ of smooth global 
sections of $E$ which reads
\begin{equation}\label{ps}
\forall (s,t)\in \Gamma(M,E)^2, \langle s,t\rangle = \int_M h_E\big(s(y), t(y)\big) |dy|.
\end{equation}
Let $P : \gme \to \gme $ be a self-adjoint elliptic  pseudo-differential operator of positive order $m$ which is bounded from below. The spectrum
of such an operator is thus real, discrete and bounded from below. Its eigenspaces
are  finite dimensional with smooth eigenfunctions, see \cite{Hormander}. We
set, for every $L\in \R$,
$$  \ul = \bigoplus_{\lambda \leq L} \ker (P-\lambda Id).$$
The dimension $N_L$ of $\ul$ satisfies Weyl's asymptotic
law
$$\frac{1}{L^{\frac{n}{m}}} N_L \toL \upn Vol \{\xi \in T^*M, \ | \ \syp (\xi)   \leq 1\},$$
where $\syp$ denotes the homogenized principal symbol of $P$,
see \cite{Hormander} and Definition $A.8$ of \cite{GaWe6}. 
The space $\ul$ inherits by restriction the $L^2$-scalar product (\ref{ps}) and 
its associated Gaussian measure defined by the density 
\begin{equation}\label{gauss}
\forall s\in \ul, \ d\mu(s)= \frac{1}{\sqrt \pi^{N_L}}
\exp( -\|s\|^2) |ds|,
\end{equation}
where $|ds|$ denotes the Lebesgue measure of $\ul$ 
associated to its scalar product. 
The measure of the discriminant $$\Delta_L = \{s\in \ul, \ s \text{ does not vanish
transversally} \} $$
vanishes  when 
$L$ is large enough, see Lemma A.1 of \cite{GaWe6}.

Our purpose is to study the topology of the vanishing locus $\sz\subset M$ 
of a section $s\in \ul$ taken at random.
More precisely, for every closed hypersurface $\Si$ of $\R^n$
not necessarily connected, and every $s\in \uld$, 
we denote by $N_\Si (s)$ the maximal number of disjoint
open subsets of $M$ with the property that every such
open subset $U'$ contains a hypersurface $\Si'$ such that 
$\Si' \subset \sz$ and $(U',\Si')$ gets diffeomorphic
to $(\R^n , \Sigma)$ (compare \cite{GaWe4}). We then set 
\begin{equation}\label{ens}
\mathbb E (N_\Si)= \int_{\ul\setminus \dl} N_\Si(s) d\mu(s)
\end{equation}
the mathematical expectation of the function $N_\Si$. Note 
that
when $\Si$ is connected, 
 the expected number of connected components 
 diffeomorphic to $\Si$
of the vanishing locus of a random section of $\ul$
gets bounded from below by $\mathbb E(N_\Si)$.

\begin{Theo}\label{Theorem 1} Let $M$ be a  smooth closed manifold of positive dimension
$n$, equipped with a Lebesgue measure $|dy|$. Let $E$ be a real
line bundle over $M$ equipped with a Riemannian metric $h_E$.
Let $P : \gme \to \gme$ be an elliptic pseudo-differential operator of 
positive order $m$, which 
is self-adjoint and bounded from below. 
Let $\Si$ be a closed hypersurface of $\R^n$, not necessarily connected.
Then, there exists a positive constant $c_\Si(P)$, such that
$$ \liminf_{L\to +\infty} \lnmm \mathbb E (N_\Si) \geq c_\Si(P).$$
\end{Theo}
The constant $c_\Si(P)$ is in fact explicit, given  by (\ref{csp}). 

Now, as in \cite{GaWe4}, we denote by $\bh_n$ the space of  diffeomorphism classes
 of closed connected  hypersurfaces of $\rn$.
For every $[\Si]\in \bh_n$ and every $ \izn$, we denote by
$ b_i(\Si)= \dim H_i(\Si,\R)$ the 
$i-$th Betti number of $\Sigma$ with real coefficients.
 Likewise, for every $s\in \ul\setminus \dl$, 
$b_i(\sz)$ denotes the  $i-$th Betti number of $\sz$, and we set
\begin{equation}\label{ebi}
\mathbb E (b_i)= \int_{\ul\setminus \dl} b_i(\sz) d\mu(s)
\end{equation}
its mathematical expectation. 

\begin{Coro}\label{Coro 1}
Let $M$ be a smooth closed manifold of positive dimension $n$ 
equipped with a Lebesgue measure $|dy|$. Let $E$ be a real line bundle
over $M$ equipped with a Riemannian metric $h_E$. Let $P : \gme \to \gme$ be
an elliptic pseudo-differential operator of positive order $m$, which is  self-adjoint and 
bounded from below. Then, for every $\izn$,
$$ \liminf_{L\to \infty} \lnmm \mathbb E(b_i) \geq
\sum_{[\Si]\in \bh_n}\sup_{\Si \in [\Si]} \big(c_\Si(P)\big) b_i(\Si).$$
\end{Coro}
Note that an upper estimate for $\mathbb E(b_i)$
of the same order in $L$ is given by Theorem 0.2 of \cite{GaWe6}.

Theorem \ref{Theorem 1} is in fact the consequence of 
Theorem \ref{Theorem 2}, which is local and
more precise. 
Let $Met_{|dy|} (M)$ be the space of Riemannian
metrics of $M$ whose associated Lebesgue measure equals $|dy|$.
For every $g\in \met$,  every  $R>0$ and every point $x\in M$, we set
\begin{eqnarray}\label{prox}
Prob^x_{\Si}(R) = \mu \Big\{  s\in \uld \ | \ \big(\sz\cap \brl\big) \supset \Si_L   \nonumber \\
 \text{ with } \big(\brl, \Si_L\big) \text{  diffeomorphic to } (\rn, \Si)\Big\},
\end{eqnarray}
where $\brl $
denotes the ball centered at $x$ of radius
$RL^{-\frac{1}m}$ for the metric $g$.

\begin{Theo}\label{Theorem 2}
Under the hypotheses of Theorem \ref{Theorem 1}, let $g\in \met$. Then,
for every $x\in M$ and every  $R>0$, 
$$\liminf_{L\to +\infty} Prob^x_{\Si}(R) \geq \pxs(R),$$
where for $R$ large enough, $p_\Si(R)=\inf_{x\in M} \pxs(R)$ 
is positive. 
\end{Theo}
Again, the function $p_\Si$ is explicit, defined by 
(\ref{rho}) (see also (\ref{rhowf}) and (\ref{rks})).
In particular, when 
$\Sigma $ is diffeomorphic to the product of spheres $ S^i \times S^{n-i-1}$, 
 Theorem \ref{theorem 3} 
provides
explicit lower estimates for the constants 
$c_\Si (P)$ and $p_\Si(R)$ appearing in 
Theorems \ref{Theorem 1} and \ref{Theorem 2}.
\begin{Theo}\label{theorem 3}
Under the hypotheses of Theorem \ref{Theorem 1}, let $g\in Met_{|dy|} (M)$
and $c_{P,g}>0$, $d_{P,g}>0$ such that 
for every $\xi \in T^*M$,
$$ d^{-1}_{P,g} \leq \frac{\si_P(\xi)^{\frac{1}m}}{\|\xi\|} \leq c^{-1}_{P,g}.$$
Then, for every $i\in \{0, \cdots, n-1\}$
and every $R\geq \frac{48\sqrt 5 n}{c_{P,g}}$,
\beq
c_{S^i\times S^{n-i-1}}(P)
& \geq &
\frac{e^{-(2\tau+1)^2}}
{2^{n+1} \sqrt \pi Vol(B(0,48 \sqrt 5n)) }
c_{P,g}^nVol_{|dy|} (M)\\
\text{and }  p_{S^i\times S^{n-i-1}}(R)&\geq &
\frac{1}{2\sqrt \pi} \exp \big(-(2\tau+1)^2\big),\\
\text{where }
\tau &=& {20}
\frac{(n+6)^{11/2}}{\sqrt {\Gamma(\frac{n}2+1)}}
\big(48n\frac{d_{P,g} }{c_{P,g}}\big)^{\frac{n+2}{2}}
\exp\Big(48\sqrt 5 n^{3/2} \frac{d_{P,g}}{c_{P,g}}\Big).
\eeq
\end{Theo}
In the case of  Laplace-Beltrami operators, we
get in particular the following.
\begin{Coro}\label{coro lap}
Let $(M,g)$
be a smooth closed $n$-dimensional 
Riemannian manifold and $\Delta$ be
its associated Laplace-Beltrami operator
acting on  functions. 
Then  for every $i\in \{0, \cdots, n-1\}$,
$$\liminf_{L\to +\infty} \frac{1}{\sqrt L^n }\mathbb E(b_i) 
 \geq c_{S^i \times S^{n-i-1}} (\Delta)
 \geq
\exp\big(-\exp(257n^{3/2})\big) Vol_g(M).$$
\end{Coro}
As a second example, Theorem \ref{theorem 3}  specializes to the case of the Dirichlet-to-Neumann operator
 on the boundary $M$ of some $(n+1)$-dimensional compact Riemannian manifold $(W, g)$. 
\begin{Coro}\label{coro diri}
Let $(W,g)$ be a smooth compact Riemannian manifold of
dimension $n+1$ with boundary $M$
and  let $\Lambda_g$ 
 be the associated 
 Dirichlet-to-Neumann operator on $M$. Then, 
 for every $i \in \{0,\cdots, n-1\}$,
$$\liminf_{L\to +\infty} \frac{1}{ L^n }\mathbb E(b_i) 
 \geq c_{S^i \times S^{n-i-1}} (\Lambda_g)
 \geq
\exp\big(-\exp(257n^{3/2})\big) Vol_g(M).$$
\end{Coro}

Note that the double exponential decay in Corollaries
 \ref{coro lap} and \ref{coro diri}
has to be compared with the exponential decay observed
in Proposition 0.4 of \cite{GaWe6} and with the analogous
double exponential decay already observed 
in Corollary 1.3 of \cite{GaWe4}. 

Let us mention some related works. In \cite{NazarovSodin},
F. Nazarov and M. Sodin proved the existence of an equivalent 
of order $L$ for the expected
number of components of the vanishing locus of 
random eigenfunctions  with eigenvalue $L$  of the Laplace operator on the
round 2-sphere. In \cite{LerarioLundberg},
A. Lerario and E. Lundberg
proved, for the Laplace operator on the round $n$-sphere, the existence
of a positive constant $c$ such that 
$\mathbb E(b_0) \geq c \sqrt L^{n}$ for large 
values of $L$.
We got in \cite{GaWe6}
upper estimates for 
$\limsup_{L\to +\infty} L^{-\frac{n}m} \mathbb E(b_i)$
under the same hypotheses as Corollary \ref{Coro 1},
and previously obtained similar 
upper and lower estimates for the expected Betti numbers or $N_\Si$'s 
of random real algebraic hypersurfaces of real projective
manifolds (see \cite{GaWe2}, \cite {GaWe3}, \cite {GaWe4}, \cite{GaWe5}).
In \cite{Letendre}, T. Letendre proved, 
under the hypotheses of Corollary \ref{coro lap},
the existence of an equivalent of order  $\sqrt L^{n}$ for the mean
Euler characteristics  (for odd $n$).
Let us finally mention 
the lecture 
\cite{Sodin}, where M. Sodin announces
a convergence in probability for $b_0$ under some hypotheses,
and
 \cite{SarnakWigman}, where 
P. Sarnak and I. Wigman announce a convergence in probability for  $N_{\Si}$ in the case of  Laplace-Beltrami operators.

In the first section, we introduce 
the space of Schwartz functions of $\R^n$ whose 
Fourier transforms have supports in the compact  
$$K_x = \{\xi \in T_x^*M \  | \ \syp (\xi) \leq 1\},$$
where $x\in M$ is given and  $T^*_x M$
is identified with $\rn$   via some isometry.
This space appears to be asymptotically a local
model for the space $\ul$. Indeed,
 any function $f$ in this space can be implemented 
in $\ul$, in the 
sense that 
there exists a family of sections $(s_L\in \ul)_{L\gg 1}$ 
 whose restriction to a ball of radius 
 of order $\lumm$ centered at $x$
 converges to $f$ after rescaling,
 see Corollary \ref{Coro 3}. 
 The vanishing locus of $f$ gets then implemented
 as the vanishing locus of the sections
 $s_L$ for $L$ large enough.
The second section is devoted to the proofs of Theorems \ref{Theorem 1} and \ref{Theorem 2}, and of Corollary \ref{Coro 1}.
For this purpose we follow the approach used in \cite{GaWe4} (see also \cite{GaWe5}), which
was itself partially inspired by the works  \cite{NazarovSodin} and \cite{Gayet}, see also \cite{LerarioLundberg}.
We begin by estimating
the expected local $C^1$-norm of elements of $\ul$, 
see Proposition \ref{Propo 2}, and then
 compare it with the amount of transversality of $s_L$.
 We can then prove Theorem
\ref{Theorem 2}, see \S \ref{II2}, and finally Theorem \ref{Theorem 1}
and its Corollary \ref{Coro 1}, see \S \ref{II3}.
The last section is devoted to the explicit estimates
and the proofs of Theorem \ref{theorem 3}
and Corollaries \ref{coro lap} and \ref{coro diri}.

\textit{Aknowledgements.} We are grateful to Olivier Druet for useful
discussions. The research leading to these results has received funding
from the European Community's Seventh Framework Progamme 
([FP7/2007-2013] [FP7/2007-2011]) under
grant agreement $\text{n}\textsuperscript{o}$ [258204].

\tableofcontents

\section{The local model and its implementation}

In the first paragraph of this section, we associate
to any closed hypersurface $\Sigma$ of $\rn$
and any symmetric compact subset $K$ of $ \R^n$
with the origin in its interior, 
a Schwartz function $f$ vanishing transversally 
along a hypersurface isotopic to $\Sigma$
and whose Fourier transform has 
support in $K$.
In the third paragraph,
we implement the function $f$ in the
neighbourhood of 
every point $\xz$ in $M$,
as the limit after rescaling 
of a sequence of sections of $\ul$.
Here,  $K$ is the pull-back
of $K_{\xz} $ under some
measure-preserving  isomorphism between $T^*_{\xz}M$
and $\R^n$. 
As a consequence, 
these sections of $\ul$ vanish 
in a neighbourhood $U_{\xz}$ of $\xz$
along a hypersurface $\Si_L$ 
of $M$ such that the pair $(U_{\xz}, \Si_L)
$ gets 
diffeomorphic
to $(\R^n, \Sigma)$.
The second paragraph 
quantifies the transversality of the vanishing
of the function $f$ and  thus of the associated sequence of
sections, in order to prepare
the estimates of the second section which involve
perturbations.

\subsection{The local model}\label{1.1}
Let $K$ be a measurable  subset of $\rn$ and let 
 $\ck$ be its characteristic function, so that 
$\ck(\xi) =1$ if $\xi\in K$ and $\ck(\xi)=0$ 
otherwise. It provides the
 projector $f \in \ldn \mapsto \ck f \in \ldn$.
After conjugation by the Fourier transform
$ \bef $ of $ \ldn $, defined for every $f\in \ldn$
and every $\xi\in \rn$ by 
$$\bef(f)(\xi) = \int_{\rn} e^{-i\cg y,\xi\cd } f(y) dy\in \ldn,$$
we get
the projector $\pi_K : \ldn  \to  \ldn,$
defined for every $ f\in \ldn$ and every $x\in \R^n$ by
$$
\pi_K(f)(x) =  \upn \int_{\xi \in K} \int_{y\in \rn}
e^{i\cg x-y,\xi \cd} f(y) d\xi dy.
$$
Note that for  $K=\rn$, $\pi_K$ is the identity map.
Denote by $\lkn$ the image of $\pi_K$. This is 
a Hilbert subspace of $\ldn$, the kernel 
of the continuous operator $Id - \pi_K = \pi_{\rn \setminus K}$.
Denote by $\cik$ the space of smooth functions on $\rn$
whose support is included in $K$, by $\srn$ the space 
of Schwartz functions of $\rn $ and set
\beqr \label{beta}
\skn=\befm (\cik).
\eeqr

\begin{Lemma}\label{Lemma 1}
Let $K$ be a bounded measurable subset of $\rn$. 
Then, $\skn \subset \lkn \cap \srn$.
\end{Lemma}
\bpr 
Since $K$ is bounded, $\cik \subset \srn$ so that
$\skn \subset \befm (\srn) = \srn$. Likewise,
for every $f\in \cik$, $\ck f = f$, so that by definition,
$f\in \lkn$. 
\epr

\newcommand{\tsi}{\widetilde \Si}
\newcommand{\ciz}{C^\infty_0(\rn)}
\begin{Lemma}\label{Lemma 2}
Let $\Si$ be a closed hypersurface of $\rn$, not necessarily connected,
and $K$ be a bounded measurable subset of $\rn$, 
symmetric with respect to the origin and which contains
the origin in its interior.
Then, there exists a hypersurface $\tsi$ of $\rn$, isotopic to $\Si$,
and a function $f_\Si $ in $\skn$ such that $f_\Si$ vanishes transversally
along $\tsi$.
\end{Lemma}
Recall that $\tsi$ is said to be isotopic to $\Si$ if and only if there exists
a continuous family $(\phi_t)_{t\in [0,1]}$ of diffeomorphisms
of $\rn$ such that $\phi_0 = Id$ and $\phi_1 (\Si)= \tsi.$

\bpr 
Let $f\in C^\infty_0(\rn)$ be a smooth compactly
supported  function of $\rn$ which vanishes transversally along $\Si$ 
and let $\tilde \chi \in \ciz$ be an even function which equals 1 in a neighbourhood of the origin. 
For every $R>0$, we set 
$$ \tilde \chi_R  : \xi \in \rn \mapsto \tilde \chi (\xi R^{-1})\in \R.$$
Then $\bef (f) \in \srn$ and $\tilde \chi_R \bef (f) $ converges
to $\bef (f)$ in $\srn$ as $R$ grows to  infinity.
Thus, $\befm (\tilde \chi_R \bef (f) )$ converges to $f$ 
in $\srn$ as $R$ grows to infinity,
and $\bef^{-1}(\tilde \chi_R \bef (f))$ takes real values.
We deduce
that when $R$ is large enough, the function
$f_R = \befm (\tilde \chi_R \bef (f))
$ is real and 
vanishes transversally in a neighbourhood of $\Si$ along
a hypersurface isotopic to $\Si$. 
By construction, the support of $\bef (f_R)$ is compact.
By hypotheses, there exists thus $\rho>0$
such that the function $\bef_\rho (f_R) : 
\xi \in \R^n \mapsto \bef (f_R) (\frac{\xi}{\rho})\in \R$
has compact support in $K$. The function
$f_\Si= \bef^{-1} (\bef_\rho (f_R))$ then belongs to
$S_K(\R^n)$ and vanishes transversally along
a hypersurface isotopic to $\Si$.  
\epr
\subsection{Quantitative transversality }\label{I2}

We now proceed as in \cite{GaWe4} 
to introduce our needed quantitative transversality estimates.
\begin{Defi}\label{Defi 1}
Let $W$ be a bounded open subset of $\rn$ and $f\in \srn$, $n>0$.
The pair $(W,f)$ is said to be regular if and only if
 zero is a regular value of the restriction of $f$ to $W$ and
the vanishing locus of $f$ in $W$ is compact.
\end{Defi}
\begin{Exam}\label{Exam 1} Let $f_\Si \in \skn\subset \srn$ be a function given
by Lemma \ref{Lemma 2}. Then, there exists a tubular neighbourhood $W$
of $\tilde \Si \subset f_\Si^{-1} (0)$ such that $(W,f_\Si)$ 
is a regular pair in the sense of Definition \ref{Defi 1}.
\end{Exam}
\begin{Defi}\label{Defi 2}
For every regular pair $(W,f)$ given by Definition \ref{Defi 1}, 
we denote by $\twf$ the set of pairs $(\delta, \epsilon) \in (\Rpe)^2
$ such that 
\beu
\item there exists a compact subset $K_W$ of $W$ such that
$\inf_{W\setminus K_W} |f|>\delta$
\item $\forall z\in W, \  |f(z )| \leq \delta \Rightarrow 
\| d_{|z}f \| >\epsilon $, where $\| d_{|z}f \|^2 = \sumun |\pif|^2(z).$
\eeu
\end{Defi}

The quantities and functions that are going to appear
in the proof of our theorems are the following.
Let $K$ be a bounded measurable subset of $\rn$. We set, for every positive 
$R$ and every $j\in \un$,
\beqr
 \rkr &=& \frac{\sqrt 2 \pedn }{\sdpn}
\inf_{t\in \Rpe} \left( \Big(\frac{R+t}{t}\Big)^{\dn} \sum^{\pedn}_{i=0}
\fr{t^i}{i!} \Big(\sum_{\substack{(j_1, \cdots, j_i)\\ \in \un^i}}
\int_K \prod_{k=1}^i|\xi_{j_k}|^2 |d\xi|\Big)^\ud
\right) \label{RR} \\ 
 \tkjr &=& \frac{\sqrt 2 {\pedn} }{\sdpn}
\inf_{t\in \Rpe} \left( \Big(\frac{R+t}{t}\Big)^{\dn} \sum^{\pedn}_{i=0}
\fr{t^i}{i!} \Big(\sum_{\substack{(j_1, \cdots, j_i)\\ \in \un^i}}
\int_K  |\xi_{j}|^2 \prod_{k=1}^i|\xi_{j_k}|^2 |d\xi|\Big)^\ud
\right).\label{TR}
\eeqr
\newcommand{\nuk}{\nu(K)}
\newcommand{\dk}{d(K)}
\begin{Rema} Denoting by $\nuk = \int_K |d\xi|$ the total measure of $K$
and by $\dk = \sup_{\xi\in K} \|\xi\|$ we note that for
every $(j_1, \cdots, j_i)\in \un^i$ and every $j\in \un$,
$$ \int_K \prod_{k=1}^i|\xi_{j_k}|^2 |d\xi|
\leq
\dk^{2i} \nuk $$
and $\int_K |\xi_j|^2\prod_{k=1}^i|\xi_{j_k}|^2
\leq \dk^{2(i+1)}\nuk.$
It follows, after evaluation at $t=R$, that 
for every $j\in \un$,
\begin{eqnarray}
\rkr &\leq &\fr{1}{\spn}{\sqrt{2\nuk}\pedn} \exp\big(R\dk \sqrt n\big)\label{majrho}\label{e1.3}\\
\tkjr &\leq &\fr{1}{\spn}{\sqrt{2\nuk}\pedn \dk } \exp\big(R\dk \sqrt n\big).\label{majm}\label{e1.4}
\end{eqnarray}
\end{Rema}

For every regular pair $(W,f)$ we set 
$$\rwf = \sup_{z\in W}\| z\|$$
and for every bounded measurable subset $K$ of $\rn$,
\beqr\label{tauwf}
\tkwf &=& \|f\|_{\ldn} 
\inf_{(\delta, \epsilon)\in \mathcal T_{(W,f)}}
\Big(\frac{1}{\delta} \rho_K (R_{(W,f)})+ 
\frac{n\sqrt n}{\epsilon} \sumun \theta^j_K(\rwf)\Big)\\
\text{and } p^K_{(W,f)} &=& \frac{1}{\sqrt \pi}
\sup_{T\in [\tkwf, +\infty[} \big(1-\frac{\tkwf}{T}\big)
\int_T^{+\infty} e^{-t^2} dt \label{rhowf}.
\eeqr
\begin{Rema}\label{rem1.7} Note that
$
p^K_{(W,f)}  
\geq \frac{1}{2\sqrt \pi} \exp\big(-(2\tkwf+1)^2\big).$
\end{Rema}
Now, let $\Si$ be a closed hypersurface of $\rn$, not necessarily 
connected. 
\begin{Defi}\label{bi}
Let  $\bi^K_\Si$ be
the set of regular pairs $\wf$ given by Definition \ref{Defi 1}
such that $f\in \skn$  and 
such that the vanishing locus
of $f$ in $W$ contains a hypersurface isotopic to $\Si$ in $\rn$.
Likewise, for every $R>0$, we set
$$ \bi^{K,R}_\Si = \{\wf \in \bi^K_\Si \ | \ \rwf \leq R\}.$$
\end{Defi}
 Finally, for every positive $R$ we set 
\beqr
 \rks (R)&=& \sup_{\wf \in \bi^{K,R}_\Si}p^K_{\wf} \label{rks}.
\eeqr
\begin{Rema}\label{pos}
It follows from Lemma \ref{Lemma 2} and Example \ref{Exam 1}
that when $R$ is large enough and $K$ satisfies the hypotheses of Lemma \ref{Lemma 2}, 
$\bi^{K,R}_\Si$ is not empty, so that
$\rks (R)>0$.
\end{Rema}

\subsection{Implementation of the local model}\label{1.3}

In this paragraph, we prove
that 
for every $x_0\in M$
and 
every measure-preserving linear
isomorphism  $A$ between  $\rn$ and $T^*_{\xz} M  $,
every function $f$ in $S_{A^*K_{x_0}}(\rn)$ 
can be implemented in $\ul$ as a sequence of sections,
see Proposition \ref{Propo 1}. 
Corollary \ref{Coro 3} then estimates
the amount of transversality 
of these sections along their
vanishing locus,
in terms of the one of $f$.

\begin{Propo}\label{Propo 1}
Under the hypotheses of Corollary \ref{Coro 1}, let $x_0\in M$,
 $\phi_{x_0} : (U_{\xz}, \xz) \subset M \to (V,0)\subset \rn$ be
a measure-preserving chart and $\tilde \chi_V \in C_c^\infty(V)$
be an even function with support in $V$  which equals 1 in a neighbourhood
of $0$. Then, for every $f\in S_{(d_{|\xz} \phi_{\xz}^{-1})^* K_{\xz}} (\rn)$,
there exists a family $(s_L)_{L\in \Rpe}\in \gme$ such that
\beu
\item \label{1} for $L$ large enough, $s_L\in \ul$ and $\|s_L\|_{L^2(M)} = \|f\|_{\ldn}$
\item \label{2} the function $z\in \rn \mapsto
L^{-\frac{n}{2m}} \tilde \chi_V (L^{-\frac{1}{m}}z)(s_L \circ \phi_{\xz}^{-1})(L^{-\frac{1}{m}}z)
\in \R$ 
converges to $f$ in $\srn$.
\eeu
\end{Propo}
Note indeed that the isomorphism $(d_{|\xz} \phi_{\xz})^{-1} : \rn \to T_{\xz}M$
defines by pull-back an isomorphism 
$((d_{|\xz} \phi_{\xz})^{-1})^* : T^*_{\xz}M \to \rn$
that makes it possible 
 to identify the compact 
 \beqr 
 K_{\xz} = \{\xi \in T^*_{\xz}M
\ |  \syp (\xi)\leq 1\}
\eeqr
 with the compact 
$\big((d_{|\xz} \phi_{\xz})^{-1}\big)^*K_{\xz}
$ of $\rn$. Moreover, 
the Riemannian metric $h_E$ of $E$ given in 
the hypotheses of 	Corollary \ref{Coro 1}
provides  a trivialization of $E$ in the neighbourhood $U_{\xz}$
of $\xz$, choosing a smaller $U_{\xz}$ if necessary,
unique up to sign. This trivialization 
makes it possible  to identify
$\tilde \chi_V s_L \circ \phi^{-1}_{\xz} $ with a function from
 $V$ to $\R$.
\bpr
For every $L\in \Rpe$, we set 
$$\tilde s_L : x\in U_{\xz} \mapsto L^{\frac{n}{2m}} \tilde \chi_V (\phi_{\xz}(x)) f(L^{\frac{1}{m}} \phi_{\xz}(x))\in E_{|x}$$
that we extend by zero to a global section of $E$. 
We denote then by $s_L$ the orthogonal projection of $\tilde s_L$ in $\ul\subset L^2(M,E)$. 
This section reads
$$s_L = \cg e_L, \tilde s_L\cd = \int_M h_E\big(e_L(x,y), \tilde s_L (y)\big) |dy|,$$
where $e_L$ denotes the Schwartz kernel of the orthogonal projection onto $\ul$. Then, 
for every $z\in \rn$, $L^{-\frac{1}{m}}z $
belongs to $ V$ when $L$ is large enough and
\beq
L^{-\frac{n}{2m}} s_L \circ \phi_{\xz}^{-1} (L^{-\frac{1}{m}}z)& = &L^{-\frac{n}{2m}}
\int_M 
h_E \Big(e_L\big(\phi_{\xz}^{-1}(L^{-\frac{1}{m}}z), y\big), \tilde s_L (y) \Big)|dy|\\
& = & \int_{U_{\xz}} \tilde \chi _V \big(\phi_{\xz}(y)\big) e_L\big(\phi_{\xz}^{-1}(L^{-\frac{1}{m}}z), y\big) 
f \big(L^{\frac{1}{m}} \phi_{\xz}(y)\big)  (y)|dy|\\
&=& L^{-\frac{n}{m}} \int_{\rn} \tilde \chi_V (L^{-\frac{1}{m}}h)(\phi_{\xz}^{-1})^* e_L(L^{-\frac{1}{m}}z, L^{-\frac{1}{m}}h)f(h) |dh|,
\eeq
where we performed the substitution $h = L^{\frac{1}{m}} \phi_{x_0} (y), $ so that $ |dh| = L^{\frac{n}{m}}|dy|$. 
But from Theorem 4.4 of \cite{Hormander},  
$$L^{-\frac{n}{m}} (\phx)^* e_L( L^{-\frac{1}{m}}z, L^{-\frac{1}{m}} h)
\toL \upn \int_{K'_{\xz}} e^{i\cg z-h,\xi\cd}  |d\xi|,$$ where $K'_{\xz} = (d_{|\xz} \phx)^*K_{\xz}$. 
Moreover, there exists $\epsilon >0$ such that this convergence holds
 in
$C^{\infty} (\rn\times \rn)$ for the semi-norms family defined by
the supremum of the derivatives of the functions on the bidisc
$\bar B(\epsilon L^{\frac{1}{m}})^2$, 
where $\bar B(\epsilon L^{\frac{1}{m}})$ denotes the closed
ball of $\rn$ of radius $\epsilon L^{\frac{1}{m}}$, see
\cite{GaWe6}.
As a consequence, after perhaps taking a smaller $V$ so that
$V$ is contained in the ball of radius $\ep$, 
$$
L^{-\frac{n}{m}}\tilde \chi_V (L^{-\frac{1}{m}}h)(\phi_{\xz}^{-1})^* e_L(L^{-\frac{1}{m}}z, L^{-\frac{1}{m}}h)f(h)
\toL \upn \int_{K'_{\xz}} e^{i\cg z-h,\xi\cd} f(h) |d\xi|$$
in this same sense, 
 which implies, with $z$ fixed, a convergence
in the Schwartz space $\srn$. After integration, it follows that
$$L^{-\frac{n}{2m}} s_L \circ \phx (L^{-\frac{1}{m}}z) 
\toL\upn \int_{K'_{\xz}} e^{i\cg z,\xi\cd } \bef (f)(\xi)|d\xi|$$
in $C^\infty (\rn)$
for our family of semi-norms 
on $\bar B (\epsilon L^{\frac{1}{m}}). $
But $f\in S_{K'_{\xz}} (\R^n)$, so that 
$$\upn \int_{K'_{\xz}} e^{i\cg z,\xi\cd } \bef(f)(\xi)|d\xi| = f(z).$$
Hence, $z\mapsto L^{-\frac{n}{2m}} s_L\circ \phx(L^{-\frac{1}{m}}z)$
converges to $f$ in $\srn$, which proves the second assertion.

If $\tilde \chi_U = \tilde \chi_V \circ \phi_{\xz}$,
we deduce that $\| s_L \tilde \chi_U\|_{L^2(M)} \toL \|f\|_{\ldn}.$
We still  need to prove that $\|s_L (1-\tilde \chi_U)\|_{L^2(M)}\toL 0.$
But since $s_L$ 
is the orthogonal projection of $\tilde s_L$ onto $\ul$, 
$$\|s_L\|_{L^2(M)} \leq \|\tilde s_L\|_{L^2(M)} 
\toL \|f\|_{L^2(\rn)}.$$
The result follows.
\epr
\begin{Coro}\label{Coro 3}
Under the hypotheses of Theorem \ref{Theorem 2}, let $\xz\in M$
and $$\px : (U_{\xz}, \xz) \subset M \to (V,0)\subset \rn$$ be
a measure-preserving  chart
such that $A=d_{|x_0} \phi_{x_0}^{-1}$
is an isometry. Let $(W, f_\Si) \in \bi_\Si^{A^*K_{\xz}}$
 and $(\delta, \ep)\in \bt_{(W, f_\Si)}$, see Definitions \ref{Defi 2} and  \ref{bi}. Then, there exist $L_0\in \R$
and $(s_L)_{L\geq L_0}$ such that for every $L\geq L_0$,
\beu
\item $s_L\in \ul$ and $\| s_L\|_{L^2(M)} \toL \|f_\Si \|_{\ldn}$ 
\item \label{2)} The vanishing locus of $s_L$ contains a hypersurface $\Sigma_L$ included
in the ball $B_g (\xz, R_{\wfs}\lumm)$ such that the pair $\big(B(\xz,R_{\wfs}\lumm ), \Si_L\big)$
is diffeomorphic to the pair $(\rn, \Si)$.
\item \label{3)} There exist two neighbourhoods $K_L$ and $W_L$ of $\Sigma_L$ such that $K_L$ is compact,
$W_L$ is open, $\Si_L \subset K_L \subset W_L \subset B_g (\xz, R_{\wfs}\lumm)$,
$\inf_{W_L\setminus K_L} |s_L| > \delta \lndm$
and for every $ y\in W_L$, 
$$ |s_L(y)| < \delta \lndm \Rightarrow \|d_{|y} (s_L \circ \phx)\| > \ep L^{\fr{n+2}{2m}}.$$
\eeu
\end{Coro}
\bpr
Let $L_0\in \R $ and $(s_L)_{L\geq L_0}$ be a family given by Proposition \ref{Propo 1} for $f=f_\Si$. Then,
the first condition is satisfied and the family of functions 
$z\in B(0, R_{\wfs}) \mapsto \lndmm s_L \circ \phx (\lumm z) $ converges to $f_\Si$ 
in $C^{\infty} (B(0, R_{\wfs})). $ Let $K$ be the compact given by Definition \ref{Defi 2},
$K_L = \phx (\lumm K) $ and $W_L = \phx (\lumm W)$. The conditions \ref{2)}. and \ref{3)}.
follow from this convergence and from Definition \ref{Defi 2}.
\epr
\section{Probability of the local presence of a hypersurface}\label{II}

In this section, we follow the  method
of \cite{GaWe4} partially inspired by \cite{NazarovSodin} 
and \cite{Gayet} (see also \cite{LerarioLundberg}, \cite{GaWe5})
in order to prove Theorem \ref{Theorem 2}.
If $\Si$ is a smooth closed  hypersurface of $\rn$,
$x\in M$ 
and $s_L\in \ul$ be given by Proposition \ref{Propo 1}, vanishing transversally along
$\Si_L$
in a small ball $B(x,L^{-\frac{1}m})$, then
we decompose any random section $s\in \ul$
as $ s= as_L + \si$, where $a\in \R$ is
Gaussian and $\sigma$ 
is taken at random in the orthogonal complement of $\R s_L$ in $\ul$.
In \S \ref{II1}, we estimate the average of 
the values of $\si$ and its derivatives on  $B(x,L^{-\frac{1}m})$
see Proposition \ref{Propo 2}.
In \S \ref{II2}, we prove that with a probability at least
$p^x_\Si$ independent of $L$, $s$ vanishes in the latter ball
along a hypersurface isotopic to $\Si_L$, thanks to the quantitative estimates
of the transversality of $s_L$ given by Corollary \ref{Coro 3},
and thanks to Proposition \ref{Propo 2}.

\subsection{Expected local $C^1$-norm of sections }\label{II1}

Recall that for $x_0\in M$,
 \beqr\label{K}
 K_{\xz} = \{\xi \in T^*_{\xz}M
\ |  \syp (\xi)\leq 1\}.
\eeqr
\begin{Propo}\label{Propo 2}
Under the hypotheses of Theorem \ref{Theorem 2}, let $\xz\in M$ and 
$\px : (U_{\xz}, \xz) \subset M
\to (V,0)\subset \rn$
 be a measure-preserving map such that $A= d_{|\xz} \phx$
is an isometry. Then, for every positive $R$ and every $j\in \un$,
\beq
\limsup_{\linf} \lndmm \mathbb E \left(\|s\|_{\Linf (B_g(\xz, R\lumm)) }\right) &\leq &\rho_{A^*K_{\xz}} (R)\\
\text{and }
 \limsup_{\linf} L^{-\fr{n+2}{2m}} \mathbb E \left(\Big\|\fr{\partial (s \circ \phx)}{\partial x_j}\Big\|_{\Linf (B_g(0, R\lumm)) }\right) 
&\leq &\theta^j_{A^*K_{\xz}} (R),
\eeq
where $\rho_{A^*K_{\xz}} $ and $\theta^j_{A^*K_{\xz}} $ are defined by (\ref{RR}) and (\ref{TR}).
\end{Propo}
\bpr
Let $t\in \Rpe$. When $L$ is large enough, the ball $B(0, (R+t)\lumm)$ 
of $\rn$ gets included in $V$. From the 
Sobolev inequality (see \S 2.4 of \cite{Federer}), 
we deduce that for every $s\in \ul$, every $k>n/2$ and every $z\in B(0, R\lumm),$
$$ |s\circ \phx (z)| \leq \frac{2k}{{Vol(B(0,t\lumm))^\ud}}
\sum_{i=0}^k (t\lumm)^i \left(\frac{1}{i!} \int_{B(0,(R+t)\lumm)} |D^i(s\circ \phx)|^2 (x) |dx|\right)^{1/2},$$
where by definition, the norm of the $i-$th derivative $D^i (s\circ \phx)$
of $s\circ \phx$ satisfies
$$ i! |D^i (s\circ \phx)(x)|^2 = 
\sum_{\substack{(j_1, \cdots, j_i)\\ \in \un^i}} \Big|\fr{\partial^i}{\partial x_{j_1}\cdots \partial x_{j_i}} (s\circ \phx)(x)\Big|^2.$$
Note indeed that the metric $h_E$ of the bundle $E$ makes it possible 
to identify $s_{|U_{\xz}} $ with a real valued function 
well defined up to a sign.
As a consequence, we deduce from the Cauchy-Schwarz inequality
\beq 
\mathbb E\Big(\| s\circ \phx\|_{\Linf (B(0, R\lumm))}\Big) &\leq &
\frac{2k}{{Vol(B(0,t\lumm))^{\ud}}} 
\sum_{i=0}^k \fr{1}{i!}(t\lumm)^i  \\
&& \left( \int_{B(0,(R+t)\lumm)} i!\mathbb E (|D^i(s\circ \phx)|^2 (x)) |dx|\right)^{1/2}.
\eeq
But given $(j_1, \cdots j_i)\in \un^i$ and $z\in B(0,(R+t)\lumm)$,
we can choose an orthonormal basis $(s_1, \cdots, s_{N_L})$
of $\ul$ such that 
$\fr{\partial^i}{\partial x_{j_1}\cdots \partial x_{j_i}} (s_l\circ \phx)(z)= 0$
for every $l>1$. Since the spectral function reads
$ (x,y)\in M\times M \mapsto e_L (x,y)= \sum_{i=0}^{N_L}
s_i(x)s_i^*(y)$, we deduce, using the decomposition of $s$ in the basis
$(s_1, \cdots, s_{N_L})$,  that 
$$ \mathbb E
 \left(\Big|\fr{\partial^i}{\partial x_{j_1}\cdots \partial x_{j_i}} (s\circ \phx)\Big|^2(z) \right)
 = \left(\int_\R a^2 e^{-a^2}\frac{da}{\sqrt \pi}\right)
 \fr{\partial^{2i}}{\partial x_{j_1}\cdots \partial x_{j_i}
 \partial y_{j_1}\cdots \partial y_{j_i}} (e_L\circ \phx)(z,z).$$
 Choosing $k= \pedn$ and noting that 
$ \int_\R a^2 e^{-a^2}\frac{da}{\sqrt \pi} = \ud$,
we deduce that for $L$ large enough,
$\mathbb E\big(\| s\circ \phx\|_{\Linf (B(0, R\lumm))}\big)$ is bounded from above by
\[
\inf_{t\in \Rpe} \frac{\sqrt 2  \pedn}{{Vol(B(0,t\lumm))^\ud}} \sum_{i=0}^{\pedn} \fr{1}{i!}(t\lumm)^i  
  \Big( \int_{B(0,(R+t)\lumm)}
\sum_{\substack{(j_1, \cdots, j_i)\\ \in \un^i} }
 \fr{\partial^{2i} e_L(x,x)}{\partial x_{j_1}\cdots \partial x_{j_i}
 \partial y_{j_1}\cdots \partial y_{j_i}} |dx|\Big)^{1/2}.
\]
Likewise, for every $j\in \{1, \cdots, n\}$, 
$\mathbb E\Big(\| \frac{\partial (s\circ \phx)}{\partial z_j}\|_{\Linf (B(0, R\lumm))}\Big)$ 
gets bounded from above by 
\[
\inf_{t\in \Rpe} \frac{\sqrt 2  \pedn}{{Vol(B(0,t\lumm))^\ud}} \sum_{i=0}^{\pedn} \fr{1}{i!}(t\lumm)^i  
  \Big( \int_{B(0,(R+t)\lumm)}
\sum_{\substack{(j_1, \cdots, j_i)\\ \in \un^i} }
 \fr{\partial^{2i+2} e_L(x,x)}{\partial x_j\partial x_{j_1} \cdots \partial x_{j_i}
 \partial y_j y_{j_1}\cdots \partial y_{j_i}} |dx|\Big)^{1/2}.
\]
Now, the result is a consequence of the asymptotic
estimate 
$$ \fr{\partial^{2i} e_L(x,x)}{\partial x_{j_1}\cdots \partial x_{j_i}
 \partial y_{j_1}\cdots \partial y_{j_i}}\equiL 
 \upn L^{\fr{n+2i}{m}}
 \int_{K_0} |\xi_{j_1}|^2\cdots |\xi_{j_i}|^2 |d\xi|,$$
 see Theorem 2.3.6 of \cite{GaWe6}. We used here that
 the balls $B_g(\xz,R\lumm)$ and $\phx (B(0,R\lumm))$
 coincide at the first order in $L$.
\epr

\subsection{Proof of Theorem \ref{Theorem 2}}\label{II2}

Let $\xz\in M$, $R>0$ and $A\in Isom_g(\rn, T_{\xz}M)$. 
Let $$\px : (U_{\xz}, \xz)\subset M \to (V,0)\subset \rn$$ be
a measure-preserving map such that $A = d_{|\xz}\phx$.
Let $\wfs\in \bi_\Si^{A^*K_{\xz, R}}$,
$(\delta, \ep)\in \bt_{\wfs}$ and $(s_L)_{L\geq L_0}$ be
a family given by Corollary \ref{Coro 3} associated to $f_\Si$,
where $K_{x_0}$ is defined by (\ref{K}). Denote by $s^\perp_L$
the hyperplane orthogonal to $s_L$ in $\ul$. Then,
$$\int_{s^\perp_L} \|s\circ \phx\|_{\Linf(\bzr)} d\mu(s)\leq
\int_{\ul} \|s\circ \phx\|_{\Linf(\bzr)} d\mu(s)$$
and for every $j\in \un$,
$$\int_{s^\perp_L} 
\Big\|\fr{\partial}{\partial  x_j} (s\circ \phx)\Big\|_{\Linf(\bzr)} d\mu(s) \leq
\int_{\ul} \Big\|\fr{\partial}{\partial x_j} (s\circ \phx)\Big\|_{\Linf(\bzr)} d\mu(s),$$
compare the proof of Proposition 3.1 of \cite{GaWe4}. From Proposition
\ref{Propo 2} and Markov's inequality we deduce that 
for every $T\in \Rpe$,
$$ \mu \Big\lbrace s\in s^\perp_L \ | \ \sup_{\bgr} |s| \geq \frac{T\delta L^{\frac{n}{2m} }}
{ \|f_\Si\|_{\ldn} } \Big\rbrace 
\leq \frac{\| f_\Si\|_{\ldn}}{T\delta} 
\rho_{A^* K_{\xz}}(R_{\wfs}) + o(1)
$$
and for every $j\in \{1, \cdots, n\}$,
$$ \mu\Big\{s\in s^\perp_L \ | \ \sup_{\bzr} |\fr{\partial}{\partial  x_j}(s\circ \phx)| 
\geq \frac{T\ep L^{\frac{n+2}{2m} }}
{ \sqrt n\|f_\Si\|_{\ldn} } \Big\} 
\leq 
\frac{\sqrt n \| f_\Si\|_{\ldn}}{T\ep} 
\theta^j_{A^* K_{\xz}}(R_{\wfs}) + o(1).
$$
It follows that the measure of the set
\beq
\be_{s^\perp_L} = \Big\{s\in s^\perp_L \ | \ \sup_{\bgr} |s|
&<& \frac{T\delta L^{\frac{n}{2m} }}
{ \|f_\Si\|_{\ldn} } \text { and } \\
 \sup_{\bzr} |d(s\circ \phx)|&<&
\frac{T\ep L^{\frac{n+2}{2m} }}
{ \|f_\Si\|_{\ldn} } \Big\}
\eeq
satisfies
$$\mu(\be_{s^\perp_L}) \geq
1 - \frac{\|f_\Si\|_{\ldn}} {T}
\Big(
\frac{1}{\delta} \rho_{A^* K_{\xz}}(R_{\wfs})
+ \frac{n\sqrt n }{\ep} \sum_{j=1}^n \theta^j_{A^* K_{\xz}}(R_{\wfs})
\Big) + o(1),
$$
where the $o(1)$ term can be chosen 
independently of $x_0$ since $M$ is compact.
Taking the supremum over the pairs $(\delta, \ep)\in 
\bt_{\wfs}$ and passing to the liminf, we deduce from
(\ref{tauwf}) the estimate 
$$ \liminf_{\linf} \mu(\be_{s^\perp_L}) \geq
1-\frac{\tau_{\wfs}^{A^*{K_{\xz}}}}{T}.$$
Now, let $$\bef_T = \Big\{ a\frac{s_L}{\|s_L\|_{L^2(M)}} + 
\si \ | \ a>T \text{ and } \si \in \be_{s^\perp_L} \Big\}.$$
From Lemma 3.6 of \cite{GaWe4}, every section
$s\in \bef_T$ 
vanishes transversally in $\bgr$ along a hypersuface
$\Si_L$ such that $(\bgr, \Si_L)$ is diffeomorphic
to $(\rn, \Si)$. 
Moreover, since $\mu$ is a product measure,
$$ \liminf_{\linf} \mu(\bef_T) \geq 
\Big(\frac{1}{\sqrt \pi} \int_T^{+\infty} e^{-t^2} dt\Big)\Big(1-\frac{\tau_{\wfs}^{A^*{K_{\xz}}}}{T}\Big).$$
Taking the supremum over $T\in [\tau_{\wfs}, +\infty[$,
we deduce from (\ref{rhowf}) that
$$ \liminf_{\linf} Prob_{\xz, \Si} (R_{\wfs}) 
\geq
\liminf_{\linf} \mu(\bef_T) \geq p^{A^*{K_{\xz}}}_{\wfs}.$$
Taking the supremum over all pairs $\wfs \in \bi^{A^*K_{x_0},R}_\Si$,
see (\ref{rks}), and then over every $A\in Isom_g(\rn, T_{\xz} M)$, 
we obtain Theorem \ref{Theorem 2} by choosing 
\bq \label{rho}\pxs (R)= \sup_{A\in \iso} (p^{A^*\kx}_\Si(R)).
\eq
Indeed, from Remark \ref{pos}, this function is positive for $R$ large enough. $\square$

\subsection{Proofs of Theorem \ref{Theorem 1} and 
Corollary \ref{Coro 1}}\label{II3}
{\bf Proof of Theorem \ref{Theorem 1}}. 
Let $g\in \met$. For every point $x$ in $M$, the
supremum 
$\sup_{R\in \Rpe} \big(\frac{1}{Vol_{eucl}(B(0,R))} p^x_\Si (R)\big)
$
is achieved and we denote by $R_m(x)$ the smallest
positive real number where it is reached. 
Denote by $\tilde g$ the normalized metric  $g/R^2_m$.
For every $L$ large enough, let $\Lambda_L$ be
a maximal subset of $M$ such that the distance between any two distinct points
of $\Lambda_L$ is larger than $2\lumm$ for $\tilde g$.
The $\tilde g-$balls centered at points of $\Lambda_L$ 
and of radius $\lumm$ are disjoint, whereas the ones
of radius $2\lumm$ cover $M$. For every $s\in \uld$
and every $x\in \Lambda_L$, we set 
$ N_{x,\Si} (s) = 1$ if $B_{\tilde g}(x,\lumm)$
contains a hypersurface $\tilde \Si$ such that 
$\tilde \Si \subset \sz$ and $(B_{\tilde g}(x,\lumm), \tilde \Si)
$
is diffeomorphic to $(\rn, \Si)$, and 
$N_{x,\Si} = 0$ otherwise. Note that 
$$ \int_{\uld} N_{x,\Si}(s) d\mu(s) \equiL Prob_{x,\Si} (R).$$
Thus,
\beq
\liminf_{\linf} \lnmm \mathbb E (N_\Si) & \geq & 
\liminf_{\linf} \lnmm \int_{\uld} \Big(\sum_{x\in \Lambda_L} N_{x,\Si} (s)\Big) d\mu(s)\\
& = & \liminf_{\linf} \lnmm \sum_{x\in \Lambda_L} Prob_{x,\Si} (R_m(x))\\
& \geq & \fr{1}{2^n} \liminf_{\linf} \sum_{x\in \Lambda_L}
Vol(B_{\tilde g} (x,2\lumm)) R_m^n(x) \Big(\fr{p^x_\Si(R_m(x))}{Vol_{eucl} B(0,R_m(x))}\Big)
\eeq
by Theorem \ref{Theorem 2}. 
Hence,  we get 
\beq
\liminf_{\linf} \lnmm \mathbb E (N_\Si) & \geq & 
\fr{1}{2^n} \int_M  \sup_{R>0} \Big(\frac{p^x_\Si (R)}{Vol_{eucl}(B(0,R)) }\Big)  R_m^n (x) |dvol_{\tilde g}(x)|\\
& = & \fr{1}{2^n} \int_M  \sup_{R>0} \Big(\frac{p^x_\Si (R)}{Vol_{eucl}(B(0,R)) } \Big) |dx|.
\eeq
Theorem \ref{Theorem 1} can be deduced after taking 
the supremum over $g\in \met$ and choosing 
the quantity $\csp$ to be equal to 
\bq \label{csp} 
\csp = \udn \sup_{g\in \met} \int_M \sup_{R>0} \Big(\frac{p^x_\Si (R)}{Vol_{eucl}(B(0,R)) }\Big) |dx|.
\eq
$\square$

{\bf Proof of Corollary \ref{Coro 1}}. 
For every $i\in \zn$ and every large enough $L>0$,
\beq
\mathbb E(b_i) & = & \int_{\uld } b_i(\sz) d\mu(s)\\
& \geq & \int_{\uld} \Big(\sum_{[\Si]\in \bh_n} N_\Si(s) b_i(\Si) \Big) d\mu(s)\\
& \geq &  \sum_{[\Si]\in \bh_n}  b_i(\Si) \mathbb E (N_\Si).
\eeq
The result is a consequence of Theorem \ref{Theorem 1} 
after passing to the liminf in the latter bound.$\square$

\section{Explicit estimates}

The goal of this section is to obtain
explicit lower bounds for the
constants  $c_\Sigma (P)$
and $\inf_{x\in M} $ $ p_\Sigma^x(R)$ 
appearing in 
Theorems \ref{Theorem 1}
and  \ref{Theorem 2}, 
when  $\Si $ is diffeomorphic to the product of spheres $ S^{i+1}\times S^{n-i-1}$
(whose $i$-th Betti number is at least one).
In the first paragraph, we 
approximate quantitatively the product of 
a  polynomial function and a Gaussian one by
a function whose Fourier transform gets compact support.
We then apply this result to a particular degree four polynomial vanishing along a product of spheres 
to finally get Theorem  \ref{theorem 3}, 
Corollary \ref{coro lap} and \ref{coro diri}.

\subsection{Key estimates for the approximation } 

Let $\tilde \chi_c : \R^n \to [0,1]$ be a smooth
function with support in the ball of radius $c>0$,
such that $\chi_c = 1$ on the ball of radius $c/2$. 
For every $Q\in \R[X_1, \cdots, X_n]$ and every $\eta>0$,
we set 
\newcommand{\qce}{q^c_{\eta}}
\newcommand{\bff}{\mathcal F}
\newcommand{\G}{e^{-\frac{\|x\|^2}{2}}}
\newcommand{\tcc}{\tilde \chi_c}
\beqr\label{refer}
q : x\in \R^n &\mapsto & q (x) =Q\G\in \R \text{ and } \\
  \qce: x\in \R^n &\mapsto &\qce (x) = \upn\int_{\rn} \tilde \chi_c (\eta \xi) \bff (q(x))(\xi) e^{i\cg x,\xi \cd} |d\xi|.
  \eeqr
  Note that $\qce\in S_{B(0,c/\eta)}(\rn)$, see (\ref{beta}).

\begin{Propo}\label{Prop 3}
Let $Q = \sum_{I\in  \Nn^n} a_I x^I
\in \R[X_1, \cdots, X_n]$ 
and $c, \eta >0$. Then,
\beu
\item 
$$\| \qce - q\|_{L^\infty (\rn)} \leq
\sqrt{\lfloor n/2+1\rfloor } \big(\frac{c}{2\eta}\big)^{\frac{n-2}{2}}
e^{-\frac{1}4 (\frac{c}{2\eta})^2} \big(\sum_{I\in \Nn^n } |a_I| \sqrt{I!}\big).$$
\item For every $k\in \{1, \cdots, n\},$
$$\Big\| \frac{\partial \qce}{\partial x_k} - \frac{\partial q}{\partial x_k} 
\Big\|_{L^\infty(\rn)}
\leq  \sqrt{\lfloor n/2+3\rfloor } \big(\frac{c}{2\eta}\big)^{\frac{n}{2}}
e^{-\frac{1}4 (\frac{c}{2\eta})^2} \big(\sum_{I\in \Nn^n } |a_I| \sqrt{I!}\big).$$
\item $$\big\| \qce - q\big\|^2_{L^2(\rn)}
\leq \sqrt{2\pi}^n N(Q) \big(\sum_{I\in \Nn^n} a_I^2 I! \big) 
e^{-\frac{1}2 (\frac{c}{2\eta})^2},$$
\eeu
where $N(Q)$ denotes the number of monomials of $Q$.
\end{Propo}
\bpr For every $x\in \rn$, we have
\beq
| \qce (x) - q(x)| 
&\leq &\upn \int_{\| \xi\| \geq \frac{c}{2\eta}}
| \bff (Q\G)| (\xi) |d\xi| \\
& \leq & \upn \sum_{I\in  \Nn^n} |a_I| \int_{\| \xi\| \geq \frac{c}{2\eta}}
| \bff(x_I \G)|(\xi) |d\xi|.
\eeq
However, 
\beqr \label{bff}
\bff (x_I \G) &=& i^{|I|} \frac{\partial}{\partial \xi_I} \big(\bff
(\G)\big)  \\
& = & \sqrt{2\pi}^n i^{|I|} \frac{\partial}{\partial \xi_I}(e^{-\frac{\|\xi\|^2}{2} })\\
& = & \sqrt{2\pi}^n i^{|I|} \prod_{j=1}^n \Big(H_{i_j} (\xi_j)e^{-\frac{\xi_j^2}{2} }\Big),
 \eeqr
 where we have set $I=(i_1, \cdots, i_n)$ and $H_j$
 the $j-$th Hermite polynomial.
 We deduce from Cauchy-Schwarz inequality that
\beq
|\qce (x) - q(x)| & \leq & \frac{1}{\sqrt{2\pi}^n}
\sum_{I\in \Nn^n} |a_I| 
\Big(
\prod_{j=1}^n \int_\R H_{i_j}^2 (\xi_j) e^{-\frac{\xi_j^2}{2}}d\xi_j
\Big)^{1/2}
\Big( \int_{\|\xi\|\geq \frac{c}{2\eta}}
e^{-\frac{\|\xi\|^2}{2}} d\xi
\Big)^{1/2}\\
& \leq & 
 \frac{1}{(2\pi)^{n/4}}\big(\sum_{I\in  \Nn^n} 
 |a_I| \sqrt{I!} \big)\sqrt{Vol(S^{n-1})}
 \Big(
 \int_{\frac{c}{2\eta}}^{+\infty}
 r^{n-1} e^{-\frac{r^2}2} dr
 \Big)^{1/2},
\eeq
since for every $k\in \Nn$, the dominating coefficient 
of $H_k(\xi)$ equals $(-1)^k$, so that an integration by parts
leads to
\beqr \label{Hermite}
 \int_\R \xi^k H_k(\xi) e^{-\frac{\xi^2}2} d\xi = (-1)^k k! \sqrt{2\pi}
\eeqr
since Hermite polynomials are orthogonal to each other.
Likewise,
after integration by parts we obtain
\beq
\int_{\cep}^{+\infty} r^{n-1} e^{-\frac{r^2}2} dr
& = & [-r^{n-2}  e^{-\frac{r^2}2}]^{+\infty}_{\cep}
+ (n-2)\int^{+\infty}_{\cep} r^{n-3} e^{-\frac{r^2}2}dr \\
& \leq & \big(\cep\big)^{n-2} e^{-\frac{1}2(\cep)^2}
+ (n-2) \big(\cep\big)^{n-4} e^{-\ud (\cep)^2} + \cdots
\eeq
From the latter we deduce, when $|\cep|\geq 1$,
\beq \int_{\cep}^{+\infty} r^{n-1} e^{-\frac{r^2}2} 
& = &\lfloor \frac{n}{2} + 1\rfloor \big(\cep\big)^{n-2} e^{-\ud (\cep)^2}(n-2)(n-4) \cdots 
\eeq
Recall that
$$Vol(S^{n-1}) = \left\lbrace \begin{array}{cc}
\frac{\sqrt {2\pi}^n }{(n-2)(n-4) \cdots 2} & \text{ if } n \text{ is even}\\
\frac{\sqrt 2\sqrt {2\pi}^n }{\sqrt \pi (n-2)(n-4) \cdots 3\times 1} & \text{ if } n \text{ is odd} 
\end{array}\right.
$$
We thus finally get
$$ \| \qce - q\|_{L^\infty} 
\leq 
\sqrt{ \lfloor \frac{n}2 + 1\rfloor } \big(\cep\big)^{\frac{n-2}2} 
e^{-\frac{1}4 (\cep)^2} \big(\sum_{I\in \Nn^n} 
|a_I| \sqrt {I!}\big).$$
Likewise, for every $k\in \{1, \cdots, n\}$,
\beq
\Big\| \frac{\partial \qce}{\partial x_k} -  \frac{\partial q}{\partial x_k}\Big\|_{L^\infty(\rn)} 
& \leq &
\frac{1}{(2\pi)^n} \int_{\|\xi\|\geq \cep } |\xi_k| |\bff (Q\G)|(\xi) |d\xi|\\
& \leq & 
\frac{1}{\sqrt {2\pi}^n}\sum_{I\in \Nn^n} |a_I| 
\big( \prod_{j=1}^n \int_\R H_{i_j}^2(\xi) e^{-\frac{\xi_j^2}2}
dx_j\big)^{\ud} \big(\int_{\|\xi\| \geq \cep } 
\| \xi\|^2 e^{-\frac{\|\xi\|^2}2} d\xi\big)^{\ud}\\
& \leq & \frac{1}{(2\pi)^{n/4}}
\sum_{I\in \Nn^n} |a_I| \sqrt {I!} Vol(S^{n-1})^{\ud}
\big(\int_{\cep}^{+\infty} r^{n+1} e^{-r^2/2} dr\big)^{\ud}\\
& \leq & \sqrt{\lfloor \frac{n}2+3\rfloor} \big(\cep\big)^{n/2}
e^{-\frac{1}4 (\cep)^2} \big(\sum_{I\in \Nn^n} |a_I| \sqrt{ I!}\big).
\eeq
Lastly,
\beq
\| \qce - q\|^2_{L^2(\rn)}
&\leq &
\big\|\bff^{-1} \big(\bff (Q\G) (1-\tcc(\eta \xi))\big)\big\|^2_{L^2(\rn)}\\
& \leq & \frac{1}{(2\pi)^n}
\int_{\xi \geq \cep} | \bff (Q \G)|^2 |d\xi| 
\text{ from Plancherel's equality}\\
& =& \int_{\xi \geq \cep } 
\big| \sum_{I\in \Nn^n} i^{|I|} a_I \prod^n_{j=1}
H_{i_j} (\xi_j) e^{-\frac{\xi_j^2}2}\big|^2 |d\xi |
\text{ from } (\ref{bff})\\
& \leq & N(Q) e^{-\frac{1}2 (\cep)^2} 
\sum_{I\in  \Nn^n}
a^2_I \prod^n_{j=1} \int_\R H_{i_j}^2(\xi_j) e^{-\frac{\xi_j^2}2}
d\xi_j \text{ from Cauchy-Schwarz}\\
& \leq & \sqrt{2\pi}^n N(Q) \big(\sum_{I\in  \Nn^n}
a^2_I I! \big) e^{-\frac{1}2(\cep)^2} \text { from } (\ref{Hermite}).
\eeq
\epr

\subsection{The product of spheres}

For every $n>0$ and every $i\in \{0, \cdots, n-1\}$,
let 
\beqr \label{Qi}
Q_i : \R^{i+1} \times  \R^{n-i-1} &\to & \R \\
(x,y) & \mapsto & (\|x\|^2 - 2)^2 + \|y\|^2 -1.
\eeqr
We recall that this polynomial vanishes
in the ball of radius $\sqrt 5$
along a hypersurface diffeomorphic 
to the product of spheres 
$S^i \times S^{n-i-1} $,
see \S 2.3.2 of \cite{GaWe4}. 
Let 
$$q_i : (x,y) \in \R^{i+1}\times
\R^{n-i-1} \mapsto Q_i(x,y) e^{-\frac{1}{2}(\|x\|^2 + \|y\|^2)}\in \R.$$
This function belongs
to the Schwartz space and
has the same vanishing locus as
$Q_i$. Let us quantify the transversality of this 
vanishing.
We set
$$W= \{ (x,y) \in \R^{i+1}\times \R^{n-i-1} , \  
\|x\|^2 + \|y\|^2 \leq 5\}.$$
\begin{Lemma}\label{lemme 5}
For every $\delta \leq 1/2$, 
$$\Big(\delta e^{-5/2}, \frac{e^{-5/2}}{2} (2 - \delta)\Big)
\in \mathcal T_{(W,q_i)},$$
see Definition \ref{Defi 2}. 
\end{Lemma}
\bpr
Let $(x,y)\in \R^{i-1}\times
\R^{n-i-1} $ be such that $\|x\|^2 + \|y\|^2 \leq 5$
and $\delta \leq 1/2$. Then
\beq 
|q_i (x,y)| < \delta e^{-5/2} & \Rightarrow &
|Q_i(x,y)|< \delta \\
& \Leftrightarrow & 1-\delta < (\|x\|^2 -2)^2 + \|y\|^2 < 1+\delta\\
& \Rightarrow & 
\left\lbrace \begin{array}{l}
\|x\|^2 > 2 - \sqrt{1+\delta} >1/2 \\
\|x\|^2 -2 >1/2 \text{ or } \|y\|^2 >1/4 \text{ since } \delta \leq 1/2.
\end{array} \right.
\eeq
Moreover, for every $ j\in \{1, \cdots, i+1\}$,
\beq
\big|\frac{\partial q_i}{\partial x_j}\big|  & \geq & 
\big|\frac{\partial Q_i}{\partial x_j}\big|e^{-5/2} - |x_j| \delta
e^{-5/2}\\
& \geq & 4|x_j| \big| \|x\|^2 - 2\big| e^{-5/2} -|x_j|\delta
e^{-5/2}\\
& \geq & |x_j| e^{-5/2}\big(4 \big| \|x\|^2 - 2\big|  -\delta\big)
\eeq
and for every $ k\in \{1, \cdots, n-i-1\}$, 
\beq
\big|\frac{\partial q_i}{\partial y_k}\big|  & \geq & 
\big|\frac{\partial Q_i}{\partial y_k}\big|e^{-5/2} - |y_k| \delta
e^{-5/2}\\
& \geq & |y_k|e^{-5/2}(2 - \delta).
\eeq
Summing up, we deduce
\beq
| d_{|(x,y)} q_i|^2 & \geq & 
\|x\|^2e^{-5} \big(4\big| \|x\|^2 - 2\big| - \delta \big)^2 + 
\|y \|^2 e^{-5}(2 - \delta)^2\\
& \geq & \frac{e^{-5}}{2} \big(4\big| \|x\|^2 - 2\big|- \delta
\big)^2 + 
\|y\|^2 e^{-5}(2 - \delta)^2\\
& \geq & \frac{e^{-5}}{4}(2 - \delta)^2.
\eeq
Since on the boundary of the ball
$W$, either $\|x\|^2\geq 7/2$ or 
$\|y\|^2 \geq 3/2$, 
the values of the function $q_i$ 
are greater than $\frac{1}2 e^{-5/2}$
and
we get the result.
\epr

We now estimate  the $L^2$-norm of $q_i$.
\begin{Lemma}\label{lemme 6}
For every $i\in \{0, \cdots, n-1\}$, 
$$\|q_i\|_{L^2(\R^n) } \leq \sqrt{\frac{3}{2}}\pi^{n/4} (n+6)^2.$$
\end{Lemma}
\bpr
We have 
\beq
\|q_i\|^2_{L^2(\R^n)} &=& \int_{\R^{i+1}\times \R^{n-i-1}}
\Big( (\|x\|^2 - 2)^2 + \|y\|^2-1\big)^2 e^{-\|x\|^2 - \|y\|^2}
dx dy\\
&=& 
\int_{\R^{i+1}\times \R^{n-i-1}} \big( \|x\|^4
- 4\|x\|^2 + 3 + \|y\|^2\big)^2 e^{-\|x\|^2 - \|y\|^2} dx dy\\
& \leq & \sqrt \pi^{n-i-1}\int_{\R^{i+1}}
(\|x\|^8 + 16 \|x\|^4) e^{-\|x\|^2}dx\\
&&+ \sqrt \pi ^{i+1} \int_{\R^{n-i-1} }
\big(\|y\|^4 + 6 \|y\|^2 + 9\big) e^{-\|y\|^2} dy\\
&& + 2\left(\int_{\R^{i+1} } \|x\|^4 e^{-\|x\|^2}dx\right)
\left(\int_{\R^{n-i-1} } (\|y\|^2 +3) e^{-\|y\|^2}dy\right).
\eeq
Now, 
\beq
\int_{\R^{i+1}} 
(\|x\|^8 + 16 \|x \|^4) e^{-\|x\|^2} dx &=& \frac{1}{2}
Vol(S^i) \int_0^{+\infty} 
(t^4 + 16 t^2) t^{\frac{i-1}{2}} e^{-t} dt \\
& = & \frac{1}{2}
Vol(S^i) \Big(\Gamma (\frac{i+9}{2}) + 16 \Gamma(\frac{i+5}2)\Big)\\
& \leq & \frac{17}{2} Vol(S^i) \Gamma(\frac{i+9}{2})
\eeq
and
\beq
\int_{\R^{n-i-1} }
\big(\|y\|^4 + 6 \|y\|^2 + 9\big) e^{-\|y\|^2} dy
&=&
\ud Vol(S^{n-i-2}) \int_0^{+\infty} 
(t^2+6t+9)t^{\ud (n-i-3)} e^{-t} dt\\
& = & 
\ud Vol(S^{n-i-2})\Big( \Gamma (\frac{n-i+3}{2})+
6\Gamma(\frac{n-i+1}2) + 9 \Gamma(\frac{n-i-1}{2})\Big) \\
& \leq & \frac{25}{2} Vol(S^{n-i-2}).
\eeq
Likewise
\beq
\int_{\R^{i+1}} \|x\|^4 e^{-\|x\|^2} dx & = &
\ud Vol(S^i) \int_0^{+\infty} t^{\frac{i+3}{2}} e^{-t} dt \\
& = & \ud Vol(S^i) \Gamma(\frac{i+5}2),
\eeq
and
\beq 
\int_{\R^{n-i-1}} (\|y\|^2 + 3)e^{-\|y\|^2} dy
&= & \ud Vol(S^{n-i-2}) \int_0^{+\infty} 
(t+3) t^{\frac{n-i-3}2} e^{-t} dt\\
& = & \ud Vol(S^{n-i-2}) \big(\Gamma (\frac{n-i+1}2) + 
3\Gamma(\frac{n-i-1}2)\big)\\
& \leq & \frac{7}2 Vol(S^{n-i-2}) \Gamma (\frac{n-i+1}2).
\eeq
Finally, since 
$$Vol(S^i) = \frac{2\pi^{\frac{i+1}{2}}}{\Gamma(\frac{i+1}{2})}
\text{ and }
Vol(S^{n-i-2}) = \frac{2\pi^{\frac{n-i-1}{2}}}{\Gamma(\frac{n-i-1}{2})}, 
$$
we get 
\beq 
\|q_i \|^2_{L^2(\R^n)} & \leq &
\sqrt \pi^n \left(
17 \frac{\Gamma(\frac{i+9}2)}{\Gamma(\frac{i+1}2)}
+ 25 
\frac{\Gamma(\frac{n-i+3}2)}{\Gamma(\frac{n-i-1}2)}
+14
\frac{\Gamma(\frac{n-i+1}2)\Gamma(\frac{i+5}2)}
{\Gamma(\frac{n-i-1}2)\Gamma(\frac{i+1}2)}
\right)\\
& \leq & 
\sqrt \pi^n \Big(
\frac{17}{16} (i+7)^4 + \frac{25}4 (n-i+1)^2
+ \frac{7}4 (n-i-1)(i+3)^2\Big)\\
& \leq & \frac{3}2 \sqrt \pi^n (n+6)^4 ,
\eeq
since $n+6\geq 7$, so that 
$\frac{25}4 (n-i+1)^2 \leq \frac{25}{4\times 49} (n+6)^4$
and $\frac{7}4 (n-i-1)(i+3)^2\leq \frac{1}4 (n+6)^4.$
\epr

We now approximate $q_i$ by a function
whose Fourier transform has compact support.
For every $i\in \{0, \cdots , n-1\}$ and $c>0$,
we set 
\beqr 
q_{i,c} : x\in \R \mapsto q_{i,c}(x) &=& q_{i,\eta}^c(\eta x)\\
& = & \frac{1}{\eta^n} \int_{\R^n}
\tcc (\xi) \bff (Q_i \G) \big(\frac{\xi}\eta\big) e^{i\cg x,\xi\cd}
|d\xi| \in \R, \label{gic}
\eeqr
see (\ref{refer}).
By construction, $q_{i,c}$ belongs to the Schwartz 
space of $\rn$ and its Fourier transform
has support in the ball of radius $c$, so that 
with the notations of \S \ref{1.1},
$q_{i,c} \in S_{B(0,c)} (\R^n ).$
\begin{Coro}\label{coro 4}
For every $i\in \{0, \cdots, n-1\}$, every $c>0$ 
and every $\eta \leq \frac{c}{48n}$,
$q_{i,c}$ vanishes in the ball
$W_\eta = \{x\in \rn, \|x\|^2\leq {5}/{\eta^2} \}
$
along a hypersurface diffeomorphic to $S^i \times
S^{n-i-1}$. 
Moreover, 
$$\Big(\frac{e^{-5/2}}4, \frac{\eta}{\sqrt 2} e^{-5/2}\Big)
\in \mathcal T_{(W_\eta, q_{i,c})}$$
and $$\|q_{i,c}\|_{L^2(\R^n)} \leq \frac{3}{2\eta^{n/2}}  \pi^{n/4} (n+6)^2.$$
 \end{Coro}
 
 \bpr
 The polynomial $Q_i$ reads
 \beq
 Q_i(x,y) = \sum_{k=1}^{i+1}x^4_k + 2\sum_{1\leq j<k\leq n}
 x_j^2 x_k^2 - 4\sum_{k=1}^{i+1} x_k^2 + \sum_{k=1}^{n-i-1} y_k^2 + 3.
 \eeq
 We deduce, with the notations of Proposition \ref{Prop 3}, 
 \beq
 \sum_{I\in \Nn^n} |a_I| \sqrt {I!} &=& (i+1)\sqrt {4!} + 4 \binom{i+1}{2} + 4\sqrt 2 (i+1) + (n-i-1)\sqrt 2 + 3
\\
& \leq & 5n + 2n^2 + 8n + 3  \leq 18 n^2
\eeq
and 
\beq
\sum_{I\in \Nn^n} a^2_I I! & = & (i+1)4! + 16 {i+1 \choose 2} + 32 (i+1) + 2(n-i-1) + 9\\
& \leq & 24 n + 8n^2 + 34 n + 9 \leq  75 n^2 ,
\eeq
whereas 
\beq
N(Q_i) & = & (i+1) + {i+1 \choose 2} + (i+1) + (n-i-1) + 1\\
 & \leq & 2n+1 + \frac{n(n-1)}2 \leq 3n^2.
 \eeq
 Noting that 
$
 \sqrt {\lfloor \frac{n}2 + 1\rfloor} \leq  \sqrt {\lfloor \frac{n}2 + 3\rfloor} \leq 2\sqrt n,$
 that $ \big(\cep\big)^{\frac{n-2}2 } \leq  \big(\cep\big)^{\frac{n}2 } $
 as soon as $\cep \geq 1$, and that 
 $$ \frac{5}2 \ln n + \frac{n}2 \ln (\cep) \leq 3n (\cep) $$
 under the same hypothesis, 
 we deduce from Proposition \ref{Prop 3} that when 
 $\eta \leq \frac{c}{48n},$
$$\| q_{i,c}(x) - q_i (\eta x) \|_{L^\infty(\rn)} \leq 36 e^{-\frac{1}8 (\cep)^2} \leq 36 e^{-72 n^2}$$
and for every $k\in \{1, \cdots, n\},$
$$\Big\| \frac{ \partial q_{i,c}}{\partial x_k}(x) - 
 \eta\frac{ \partial q_i}{\partial x_k}(\eta x) \Big\|_{L^\infty(\rn)}  \leq 36
 \eta e^{-72 n^2}.$$
From Lemma \ref{lemme 5}  follows, choosing $\delta=1/2$,
 that for every $x\in \rn$
such that $\|x\|^2  \leq 5/\eta^2 $ and every $\eta \leq \frac{c}{48n}$,
\beq 
q_{i,c}(x) \leq \frac{e^{-5/2}}4 & \Rightarrow & 
q_i (\eta x) \leq \frac{e^{-5/2}}2 \\
& \Rightarrow & | d_{| \eta x} q_i | > 3\frac{e^{-5/2}}4 \\
& \Rightarrow & | d_{|x} q_{i,c}| > \eta \frac{e^{-5/2}}{\sqrt 2},
\eeq
since 
\beq 
|d_{|x} q_{i,c}|&\geq & \eta |d_{|\eta x} q_i|
- |d_{|x} q_{i,c} - \eta d_{|\eta x} q_i| \\
& > & \eta \frac{3e^{-5/2}}{4}
- \sqrt {\sum_{k=1}^n 
\Big|\frac{ \partial q_{i,c}}{\partial x_k}(x) - 
 \eta\frac{ \partial q_i}{\partial x_k}(\eta x)\Big|^2} \\
 & \geq & \eta \big(\frac{3e^{-5/2}}{4} - 
 36 \sqrt n e^{-72n^2} \big) \\
& > & \eta \frac{e^{-5/2}}{\sqrt 2}.
\eeq
From Lemma 3.6 of \cite{GaWe4}, 
$q_{i,c}$ vanishes 
in the ball $W_\eta$ along a hypersurface diffeomorphic to
$S^i \times S^{n-i-1}$ and by definition,
$(\frac{e^{-5/2}}4, \eta\frac{ e^{-5/2}}{\sqrt 2}) \in \mathcal T_{(W_\eta,q_{i,c})}$
if $\eta \leq \frac{c}{48n}.$

Lastly, we estimate the $L^2$-norm of $q_{i,c}$. By 
Proposition \ref{Prop 3}
and the bounds given above, $\| q^c_{i,\eta} - q_i \|^2_{L^2(\rn)}\leq 
\sqrt{2\pi}^n 225
n^4 e^{-288 n^2}$, so that 
\beq
\|q_{i,c}\|_{L^2(\rn)} &=& \frac{1}{\eta^{n/2}} \|q^c_{i,\eta} \|_{L^2(\rn)} \\
& \leq & \frac{1}{\eta^{n/2}}(\|q_i\|_{L^2 (\rn)}
+ \| q^c_{i,\eta} - q_i \|_{L^2(\rn)})\\
& \leq & \frac{1}{\eta^{n/2}}\Big(\sqrt{\frac{3}2}\pi^{n/4} (n+6)^2 + (\sqrt{2\pi}^n 225
n^4 e^{-288 n^2})^{1/2}\Big)
\text{ by Lemma \ \ref{lemme 6}}
\\
& \leq & \frac{3}{2\eta^{n/2}} \pi^{n/4} (n+6)^2.
\eeq
 \epr
 \subsection{Proofs of Theorem  \ref{theorem 3}, Corollary \ref{coro lap} and Corollary \ref{coro diri}}
\bpr[ of Theorem \ref{theorem 3}]
Let us choose $c=c_{P,g}$ and $\eta =  \frac{c_{P,g}}{48n}$
see the definition (\ref{gic}) of $q_{i,c}$. 
It follows from Corollary \ref{coro 4} that 
$R\geq \frac{48\sqrt 5 n }{c_{P,g}}$,
$(W_\eta,q_{i,c}) \in 
\mathcal I^{B(0,c_{P,g}),R}_{S^i \times S^{n-1-i}}$
so that for any 
$x\in M$ and any
$A\in Isom_g (\R^n, T_xM)$,
$$(W_\eta,q_{i,c}) \in \mathcal I^{A^*K_x,R}_{S^i \times S^{n-1-i}}.$$
Indeed, $$B(0,c_{P,g} ) \subset A^*K_x \subset B(0,d_{P,g}).$$
  From Remark \ref{rem1.7},
we get that for
every 
 $x\in M$ and every $R\geq \frac{48\sqrt 5 n }{c_{P,g}}$, 
$$ p^x_{S^i\times S^{n-i-1}}(R) \geq \frac{1}{2\sqrt \pi}
\exp(-(2\tau+1)^2).$$
From (\ref{e1.3}),  (\ref{e1.4}), 
(\ref{tauwf}) and Corollary \ref{coro 4}
with $\eta = \frac{c_{P,g}}{48n}$,
using that $\nu(A^*K_x) \leq Vol(B(0,d_{P,g}))$,
we deduce 
\beq
\tau & \leq & 
\frac{3}{ 2}  \pi^{n/4} (n+6)^2 \big(\frac{48n}{c_{P,g}}\big)^{n/2}
\Big(\frac{4}{e^{-5/2}} \frac{1}{\sqrt \pi^n } \sqrt {2 Vol(B(0,d_{P,g}))} 
\lfloor \frac{n}2 + 1\rfloor
\exp \big(48 \sqrt 5 n \sqrt n \frac{d_{P,g}}{c_{P,g}}\big)\\
& + & \frac{ 48n\sqrt 2}{e^{-5/2} c_{P,g}}n \sqrt n 
\frac{n}{\sqrt \pi^n }\sqrt {2 Vol(B(0,d_{P,g}))}
\lfloor \frac{n}2 + 1\rfloor  d_{P,g}
\exp \big(48 \sqrt 5 n \sqrt n \frac{d_{P,g}}{c_{P,g}}\big)\Big)\\
& \leq & \frac{3}{4\pi^{n/4}}(n+6)^3 (48n)^{n/2} \sqrt {2 Vol(B(0,1))}
(\frac{d_{P,g}}{c_{P,g}})^{n/2} \exp \big(48 \sqrt 5 n^{3/2} \frac{d_{P,g}}{c_{P,g}}\big) \\
&&\big(4e^{5/2} + \sqrt 2 e^{5/2} n^{5/2} (48n) \frac{d_{P,g}}{c_{P,g}}\big)\\
& \leq & 20 \frac{(n+6)^{11/2}}{\sqrt{\Gamma(\frac{n}2+1)}}(48n \frac{d_{P,g}}{c_{P,g}})^{\frac{n+2}2} \exp\big(48\sqrt 5 n^{3/2} \frac{d_{P,g}}{c_{P,g}}\big).
\eeq
The estimate for 
$c_{[S^i\times S^{n-i-1}]} $ follows from 
the above estimate with $R = 48\sqrt 5 \frac{n}{c_{P,g}}$,
see (\ref{csp}).
\epr

\bpr[ of Corollary \ref{coro lap}]
If $P$ is the Laplace-Beltrami operator 
associated to a metric $g$ on $M$,
then we choose as the Lebesgue measure $|dy|$ 
on  $M$ the measure $|dvol_g|$ 
associated to $g$, so that $g\in Met_{|dy|} (M)$
and the principal symbol of $P$
equals 
$\xi \in T^*M \mapsto \|\xi\|^2 \in \R.$
Theorem \ref{theorem 3} then applies with $m=2$
and $c_{P,g} = d_{P,g} = 1$ 
and we deduce, using $\Gamma(\frac{n}2 + 1)\geq
1/2$, that
\beq 
\tau & \leq & 20 \frac{(7n)^{11/2} }{\sqrt {\Gamma(\frac{n}2 + 1)} }
(48n)^{\frac{n+2}{2}} \exp(108n^{3/2})\\
& \leq & \exp \big(
\ln(20\sqrt 2) + \frac{11}2 \ln 7 + \frac{n+2}2 \ln 48 
+ \frac{13}2 \ln n + \frac{n}2 \ln n + 108 n^{3/2}\big)\\
& \leq & \exp(18 + \frac{17}2 (n-1) + \frac{n}2 (2\sqrt n - 1)
+ 108 n^{3/2} )\\
& \leq & \exp (127 n^{3/2}).
\eeq 
Theorem \ref{theorem 3} then provides 
for every $i\in \{0, \cdots, n-1\}$, 
\beq
\big(Vol_g(M)\big)^{-1}c_{[S^i \times S^{n-i-1}]} (P)
& \geq & 
\exp \big(
-(2\tau+ 1)^2 - (n+1) \ln 2 - \frac{1}2 \ln \pi \\
&&- n \ln (48 \sqrt 5 n)
- \ln (\pi^{n/2}) + \ln (\Gamma (n/2+1))\big) \\
& \geq & \exp \big( -(2\tau+1)^2 - 3/2 - 6 \ln n - n \ln n\big)\\
& \geq & \exp \Big(-\exp(256 n^{3/2}) - \exp \big(\ln (17/2) + \ln n +
\ln (\ln n) \big) \Big)\\
& \geq & \exp(-\exp (257 n^{3/2}).
\eeq
\epr
\begin{Rema} 
Under the assumptions of Corollary \ref{coro lap}, 
we get likewise 
for $R\geq 48 \sqrt 5 n $,
\beq
\inf_{x\in M} \big(p^x_{S^i\times S^{n-i-1}} (R)\big)
& \geq & \frac{1}{2\sqrt \pi } \exp (-\exp(256n^{3/2})) \\
& \geq & \exp (-\exp(257n^{3/2})).
\eeq
\end{Rema}
\bpr[ of Corollary \ref{coro diri}]
If $P$ denotes the Dirichlet-to-Neumann operator
on $M$, then
the principal symbol of $P$
equals 
$\xi \in T^*M \mapsto \|\xi\| \in \R.$
Theorem \ref{theorem 3} then applies with $m=1$
and $c_{P,g} = d_{P,g} = 1$.
Thus, the proof is the same as the one 
of Corollary \ref{coro lap}.
\epr

\providecommand{\bysame}{\leavevmode\hbox to3em{\hrulefill}\thinspace}
\providecommand{\MR}{\relax\ifhmode\unskip\space\fi MR }
\providecommand{\MRhref}[2]{%
  \href{http://www.ams.org/mathscinet-getitem?mr=#1}{#2}
}
\providecommand{\href}[2]{#2}


Damien Gayet\\
 Univ. Grenoble Alpes, IF, F-38000 Grenoble, France\\
CNRS, IF, F-38000 Grenoble, France\\
 damien.gayet@ujf-grenoble.fr\\

Jean-Yves Welschinger\\
 Universit\'e de Lyon \\
CNRS UMR 5208 \\
Universit\'e Lyon 1 \\
Institut Camille Jordan \\
43 blvd. du 11 novembre 1918 \\
F-69622 Villeurbanne cedex, France\\
 welschinger@math.univ-lyon1.fr

\end{document}